%% file: nilpotent.tex
\documentclass[reqno]{amsart}
\usepackage{setspace}

\usepackage{epsfig,amssymb,amsmath,mathrsfs,color}

\newcommand{\R}{\mathbb{R}}
\newcommand{\N}{\mathbb{N}}
\newcommand{\Z}{\mathbb{Z}}
\newcommand{\union}{\cup}
\newcommand\conv{\operatorname{conv}}

\newcommand\intersection{\cap}
\newcommand\closure{\operatorname{cl}}
\newcommand\after{\circ}

\newcommand\facetgeo{\Lambda}
\newcommand\word{w_F}
\newcommand\evalword{\overline{w}_F}

\newtheorem{prop}{Proposition}[section]

\newtheorem{corollary}[prop]{Corollary}
\newtheorem{lemma}[prop]{Lemma}

\newtheorem{theorem}[prop]{Theorem}

\begin{document}

\title[Action of a nilpotent group on its horofunction boundary]
{The action of a nilpotent group on its horofunction boundary
has finite orbits
}
\date{\today}
\author{Cormac Walsh
}
\address{
INRIA and CMAP, Ecole Polytechnique.
Postal address: CMAP, Ecole Polytechnique, 91128 Palaiseau C\'edex, France}
\email{cormac.walsh@inria.fr}
\thanks{This work was partially supported
by the joint RFBR-CNRS grant number 05-01-02807.}




\begin{abstract}
We study the action of a nilpotent group $G$ with finite generating set $S$
on its horofunction boundary. We show that there is one finite orbit associated
to each facet of the polytope obtained by projecting $S$ into the torsion-free
component of the abelianisation of $G$.
We also prove that these are the only finite orbits of Busemann points.
To finish off, we examine in detail the Heisenberg group with its usual
generators.
\end{abstract}

\newcommand\isom{\operatorname{Isom}}

\maketitle

\section{Introduction}

Given a group acting by isometries on a space, one often finds it useful
to study the induced action on some boundary of the space,
for example the hyperbolic boundary of a Gromov hyperbolic space or the
ideal boundary of a CAT(0) space.
The hope is that the action on the boundary is simpler than the action on
the space itself.

Of course, we may endow any finitely-generated group with its word metric,
and then the group acts isometrically on itself by left translations.
One can then study the induced action on the horofunction boundary,
which exists for any metric space.

The case of $\Z^d$ with an arbitrary finite generating set was investigated by
Rieffel~\cite{rieffel_group}
in connection with his work on quantum metric spaces.
Essential to his work was that the action on the horoboundary is amenable
and that there are sufficiently many finite orbits.
In~\cite{develin_cayley}, Develin studied in greater detail the horoboundary
of $\Z^d$ and the action of $\Z^d$ on it.

In this paper, we study the case of a finitely-generated nilpotent group.
We find that, just as in the case of $\Z^d$, there are always
finite orbits, whatever the generating set.

To state our results, we need the following construction.
Recall that the quotient $G/[G,G]$ of a group $G$ by its commutator subgroup
is an abelian group, called the \emph{abelianisation}.
If $G$ is finitely generated, then so is the abelianisation, which
can therefore be written as a direct product $Z_0 \times \Z^N$,
where $Z_0$ is a finite abelian group (the torsion subgroup) and $N\in\N$.
So we can define a projection map $\phi$ assigning to every element of $G$
its equivalence class in $(G/[G,G])/Z_0$.
We consider $\Z^N$ to be embedded in $\R^N$.
Let $P:=\conv(\phi(S))$ be the convex hull of the image of
the generating set $S$ under the map $\phi$.

We will assume for simplicity that $S$ is symmetric, that is, $S=S^{-1}$,
since this makes the word length distance a metric. The results should
also be true, however, in the more general case.

A particularly interesting subset of the horofunction boundary is its
set of Busemann points. These are closely related to the geodesic paths.
In general, all geodesic paths converge to Busemann points~\cite{rieffel_group}.
Moreover, when the metric is integer valued, as in the present setting,
the Busemann points are precisely the limits of the geodesic paths.

Recall that a facet of a polytope is a proper face of maximal dimension,
in other words, of co-dimension $1$.

\begin{theorem}
\label{thm:maintheorem}
Let $G$ be a nilpotent group with finite symmetric generating set $S$
and consider the action of $G$ on its horofunction boundary coming from
the word length metric associated to $S$.
Then, there is a natural one-to-one correspondence between the finite orbits
of Busemann points and the facets of $P$.
\end{theorem}

When $G=\Z^d$, the map $\phi$ is just the identity map.
In this case, it was shown in~\cite{develin_cayley} that all boundary points
are Busemann and that there is a one-to-one correspondence between
the orbits and the proper faces of $P$, with the finite orbits corresponding
to the facets.
We will give an example later of a non-abelian nilpotent group
where this stronger correspondence breaks down.

We also have more precise results for the discrete Heisenberg group,
which has the following presentation:
\begin{align*}
H_3 := \langle a,b \mid [[a,b],a] = [[a,b],b] = e \rangle.
\end{align*}
We use the word length metric coming from the usual set of
generators $S:=\{a,b,a^{-1},b^{-1}\}$.
For this group with these generators, an exact formula for the
word length metric is known.
We show directly that there are four finite orbits of Busemann points,
each consisting of a single horofunction. These horofunctions are the limits
of the four geodesic words
$a^{\epsilon_a} b^{\epsilon_b} a^{\epsilon_a} b^{\epsilon_b} \cdots$, with
$\epsilon_a =\pm 1$ and $\epsilon_b =\pm 1$.
We also show that every other Busemann point is the limit
of a geodesic path of the form $t\to c^m b^n a^{\pm t}$ or
$t\to c^{m\pm lt} b^{\pm t} a^l$, where $c:=[a,b]$.
Using the formula for the metric, all these Busemann points can be calculated
explicitly; see Theorems~\ref{thm:corners} and~\ref{thm:heisenberg_busemen}.

\section{The Horofunction Boundary}

The horofunction boundary of a proper metric space $(X,d)$ is defined
as follows. One assigns to each point $z\in X$ the function
$\psi_z:X\to \R$,
\begin{equation*}
\psi_z(x) := d(x,z)-d(b,z),
\end{equation*}
where $b$ is some base-point.
The map $\psi:X\to C(X),\, z\mapsto \psi_z$ defines a continuous injection
of $X$ into $C(X)$,
the space of continuous real-valued functions on $X$ endowed
with the topology of uniform convergence on compact sets.
The horofunction boundary is defined to be
\begin{align*}
X(\infty):=\closure\{\psi_z\mid z\in X\}\backslash\{\psi_z\mid z\in X\},
\end{align*}
and its elements are called horofunctions.
This definition seems to be due to Gromov~\cite{gromov:hyperbolicmanifolds}.
For more information, see ~\cite{bgs}, ~\cite{rieffel_group},
and~\cite{AGW-m}.

The action of the isometry group $\isom(X,d)$ of $X$ extends continuously
to an action by homeomorphisms on the horofunction compactification.
One can check that
\begin{align*}
g\cdot \xi(x) = \xi(g^{-1}\cdot x) - \xi(g^{-1}\cdot b),
\end{align*}
for any isometry $g$, point $x$, and horofunction $\xi$.

The metric spaces in which we will be interested in this paper are
the finitely-generated nilpotent groups with their word length metrics.
Recall that the Cayley graph of a group $G$ with symmetric generating set $S$
has a vertex for each element of $G$
and an edge between two vertices $u$ and $v$ if $u^{-1}v\in S$.
The word length distance $d(u,v)$ between $u$ and $v$ in $G$ is the length
of the shortest path between the corresponding vertices in the Cayley graph.
By construction, it is a left invariant metric.
We choose the identity element $e$ as the base-point.

A path is a sequence of group elements. We say that a path is
\emph{geodesic} if it is an isometry from $\N$ to $G$.
It is known that geodesic paths converge to points in the horofunction
boundary.
In the present setting, since the metric is integer valued, we may define a
\emph{Busemann point} as the limit of a geodesic. See~\cite{rieffel_group}
and~\cite{AGW-m} for the definition in more general spaces.

For any finite word $w$ with letters in $S$, we denote by $\overline w$
the element of $G$ obtained by multiplying the letters of $w$,
and we say that $w$ \emph{represents} $\overline w$.
We denote by $|w|$ the length of $w$, that is, its number of letters.
We say that a finite word $w$ is a \emph{geodesic word}
if $|w|=d(e,\overline w)$, and that an infinite word is geodesic if each of
its finite prefixes is. Obviously, there is a one-to-one correspondence
between geodesic paths and pairs $(g,w)$, where $g$ is in $G$ (the starting
point of the path) and $w$ is an infinite geodesic word.
The \emph{steps} of a geodesic path are the letters of the corresponding
geodesic word. Continuing to use terminology of~\cite{develin_cayley},
we say that the \emph{direction} of a geodesic path is the set of letters
that occur infinitely often as steps.

We start off with an easy proposition telling us when two geodesics converge
to the same Busemann point.
Note that the proof relies on the fact that the horofunctions are
integer valued.
\begin{prop}
\label{prop:switch}
Let $\gamma_1$ and $\gamma_2$ be two geodesics in a finitely-generated group.
The following are equivalent:
\begin{itemize}
\renewcommand{\labelitemi}{(i)}
\item
$\gamma_1$ and $\gamma_2$ converge to the same Busemann point;
\renewcommand{\labelitemi}{(ii)}
\item
for all $N\in\N$, there is a geodesic agreeing with $\gamma_1$ up to time $N$
that has a point in common with $\gamma_2$ and eventually agrees with
$\gamma_1$ again;
\renewcommand{\labelitemi}{(iii)}
\item
there is a geodesic having infinitely many points in common with both
$\gamma_1$ and $\gamma_2$.
\end{itemize}
Moreover, the above statements hold if
\begin{itemize}
\renewcommand{\labelitemi}{(iv)}
\item
$\gamma_1$ and $\gamma_2$ have a point in common, and, for all $N\in\N$,
there is a geodesic agreeing with $\gamma_1$ up to time $N$ with the same limit
as $\gamma_2$, and a geodesic agreeing with $\gamma_2$ up to time $N$ with
the same limit as $\gamma_1$
\end{itemize}
\end{prop}
\begin{proof}
The implications (ii)$\Rightarrow$(iii) and (iii)$\Rightarrow$(i) are obvious.

Assume (i) holds, that is, $\gamma_1$ and $\gamma_2$ converge to some
horofunction $\xi$. Let $N\in\N$.
We can find $t_2$ large enough that
\begin{align*}
\xi(\gamma_2(0)) - \xi(\gamma_1(N))
   = d(\gamma_2(0),\gamma_2(t_2)) - d(\gamma_1(N), \gamma_2(t_2)).
\end{align*}
We can then find $t_1$ large enough that
\begin{align*}
\xi(\gamma_2(t_2)) - \xi(\gamma_1(N))
   = d(\gamma_2(t_2),\gamma_1(t_1)) - d(\gamma_1(N), \gamma_1(t_1)).
\end{align*}
Since $\gamma_2$ is a geodesic,
\begin{align*}
\xi(\gamma_2(0)) = d(\gamma_2(0),\gamma_2(t_2)) + \xi(\gamma_2(t_2)).
\end{align*}
It follows that
\begin{align*}
d(\gamma_1(N), \gamma_1(t_1))
   = d(\gamma_1(N), \gamma_2(t_2)) + d(\gamma_2(t_2),\gamma_1(t_1)).
\end{align*}
So there exists a geodesic agreeing with $\gamma_1$ up to time $N$,
passing through $\gamma_2(t_2)$, and agreeing with $\gamma_1$ after time $t_1$.
This establishes (ii).

Now assume that (iv) holds. So there is a geodesic $\gamma_3$ agreeing with
$\gamma_1$ up to any given time $N$ with the same limit as $\gamma_2$.
Applying the equivalence of (i) and (ii) just proved, we see there is a
geodesic $\gamma_4$ agreeing with $\gamma_1$ up to time $N$ and agreeing with
$\gamma_2$ after some time $M$.
A symmetrical argument shows that there is a geodesic $\gamma_5$
agreeing with $\gamma_2$ up to time $M$ and eventually agreeing with $\gamma_1$.
It is easy to prove that the path agreeing with $\gamma_4$ up to time $M$ and
thereafter agreeing with $\gamma_5$ is a geodesic satisfying the requirements
of (ii).
\end{proof}

\section{Preliminaries on Nilpotent Groups}

We use $[a,b]$ to denote the commutator $a^{-1}b^{-1}ab$ of two elements
of a group. Also, for any two subgroups $A$ and $B$, we denote by
$[A,B]$ the subgroup generated by all $[a,b]$ with $a\in A$ and $b\in B$.
For any group $G$, the subgroup $[G,G]$ is called the
\emph{commutator subgroup}. It is always normal.

Recall that a group $G$ is \emph{nilpotent} if $\Gamma_n$ is trivial for
some $n\in\N$, where $\Gamma_i$; $i\in\N$ is the \emph{lower central series}
defined by $\Gamma_1:=G$ and $\Gamma_{i+1}:= [\Gamma_i,G]$ for $i\in\N$.

We will need the following two lemmas concerning finite-index subgroups
of nilpotent groups. They have similar proofs, which rely on the identities
\begin{align}
\label{eqn:commuteright}
[x,yz] &= [x,z] [x,y] [[x,y],z] \qquad\text{and} \\
\label{eqn:commuteleft}
[xy,z] &= [x,z] [[x,z],y] [y,z].
\end{align}

\begin{lemma}
\label{lem:comm_finite_index}
Let $H$ be a finite index subgroup of a finitely generated nilpotent group~$G$.
Then, $[H,H]$ has finite index in $[G,G]$.
\end{lemma}
\begin{proof}
We use an inductive argument.
The induction hypothesis is that there is a finite set $R_i\subset G$ such that
any $g\in[G,G]$ can be expressed as $g=h'r' \mod \Gamma_i$,
with $h'\in [H,H]$ and $r'\in R_i$.

The case when $i=2$ is trivial.

Now assume the hypothesis is true for some $i\ge 2$.
The assumption of the lemma implies
that, for each generator $s\in S$, there is some $n_s\ge 1$ and $h_s\in H$
such that $s^{n_s} = h_s$.
Let $C_i$ denote the set of simple commutators of weight $i$ of elements of the
generating set $S$.
In other words, $C_1 := S$ and
\begin{align*}
C_{j+1} := \big\{ [c,s] \mid \text{$c\in C_j$ and $s\in S$} \big\},
\qquad \text{for $j\ge 1$}.
\end{align*}
Let $c := [\dots[s_1, \dots] , s_i ]$ be an arbitrary element of $C_i$.
Using the identities~(\ref{eqn:commuteright}) and~(\ref{eqn:commuteleft}),
one can show that
\begin{align*}
c^{n_{s_1}\cdots n_{s_i}}
   &= [\dots[s_1^{n_{s_1}},\dots],s_i^{n_{s_i}}] \qquad\mod \Gamma_{i+1}  \\
   &= [\dots[h_{s_1} ,\dots], h_{s_i}] \qquad\mod \Gamma_{i+1}.
\end{align*}
So, letting $p_c:= n_{s_1}\cdots n_{s_i}$,
we see that $c^{p_c}$ is in $[H,H]$, modulo $\Gamma_{i+1}$.

As $\Gamma_i/\Gamma_{i+1}$ is abelian and generated by the elements of
$C_i \mod \Gamma_{i+1}$ we can write any $d\in \Gamma_i$ in the form
\begin{align*}
d = \prod_{c\in C_i} c^{m_c}    \qquad\mod\Gamma_{i+1},
\end{align*}
where the product is taken in any fixed order and $m_c\in \Z$ for all
$c\in C_i$. Therefore, from what we have just proved,
\begin{align*}
d = h \prod_{c\in C_i} c^{q_c}    \qquad\mod\Gamma_{i+1},
\end{align*}
for some $h\in[H,H]$, where $q_c$ is the remainder when $m_c$ is divided by
$p_c$. Let $R$ be the finite set
\begin{align*}
R:= \Big\{\prod_{c\in C_i} c^{q_c}
        \mid \text{$0\le q_c < p_c$ for all $c\in C_i$} \Big\}.
\end{align*}

By the induction hypothesis, there is a finite set $R_i\subset G$ such that
any $g\in[G,G]$ can be expressed as $g=h'r'd$, with $h'\in [H,H]$, $r'\in R_i$,
and $d\in \Gamma_i$.
So $g=h'dr' \mod \Gamma_{i+1}$. By the result of the previous paragraph,
$d =hr \mod \Gamma_{i+1}$ for some $h\in [H,H]$ and $r\in R$.
Therefore, $g=h'hrr' \mod \Gamma_{i+1}$, and so the induction hypothesis
is true for $i+1$.
\end{proof}

\begin{lemma}
\label{lem:group_finite_index}
Let $H$ be a subgroup of a finitely generated nilpotent group $G$,
and let $\Theta$ be the natural homomorphism from $G$ to its abelianisation
$G/[G,G]$. If $\Theta(H)$ has finite index in $\Theta(G) = G/[G,G]$, then
$H$ has finite index in $G$.
\end{lemma}
\begin{proof}
We use again an inductive argument.
The induction hypothesis here is that there is a finite set $R'\subset G$
such that any $g\in G$ can be expressed as $g=h'r' \mod \Gamma_i$,
with $h'\in H$ and $r'\in R'$.

By assumption, this is true for $i=2$.

Now assume it is true for some $i\ge 2$. The assumption of the lemma implies
that, for each generator $s\in S$, there is some $n_s\ge 1$ and $h_s\in H$
such that $\Theta(s^{n_s}) = \Theta(h_s)$.
Therefore, for all $s\in S$, there exists $g_s\in[G,G]$ such that
$s^{n_s} = h_s g_s$.

As in the proof of the previous lemma, let $C_i$ denote the set of simple
commutators of weight $i$ of elements of $S$, and consider an arbitrary
element $c := [\dots[s_1, \dots] , s_i ]$.
Using the identities~(\ref{eqn:commuteright}) and~(\ref{eqn:commuteleft}),
one can show as before that
\begin{align*}
c^{p_c}
   &= [\dots[s_1^{n_{s_1}} ,\dots],s_i^{n_{s_i}}] \qquad\mod \Gamma_{i+1}  \\
   &= [\dots[h_{s_1} ,\dots], h_{s_i}] \qquad\mod \Gamma_{i+1},
\end{align*}
where $p_c:= n_{s_1}\cdots n_{s_i}$.

It follows just as before that any $d\in \Gamma_i$ can be written as
$d =hr \mod \Gamma_{i+1}$ for some $h\in H$ and $r\in R$,
where $R$ is the finite set
\begin{align*}
R:= \Big\{\prod_{c\in C_i} c^{q_c}
        \mid \text{$0\le q_c < p_c$ for all $c\in C_i$} \Big\}.
\end{align*}

By the induction hypothesis there is a finite set $R'\subset G$ such that
any $g\in G$ can be expressed as $g=h'r'd$, with $h'\in H$, $r'\in R'$, and
$d\in \Gamma_i$.
So $g=h'dr' \mod \Gamma_{i+1}$.
But we have seen that $d =hr \mod \Gamma_{i+1}$ for some $h\in H$ and $r\in R$,
and therefore, $g=h'hrr' \mod \Gamma_{i+1}$.
This completes the induction step.
\end{proof}

\section{Facets and Finite Orbits}

Let $G$ be a nilpotent group generated by a finite symmetric set $S$.
Recall that $\phi$ is the projection onto the torsion-free component of the
abelianisation of $G$, and that $P$ is the convex hull of $\phi(S)$,
considered as a subset of $\R^N$.

\begin{lemma}
\label{lem:geo_facet}
Let $V$ be a subset of $S$ such that $\phi(V)$ is contained
in a facet of $P$. Then, any word with letters in $V$ is geodesic.
\end{lemma}
\begin{proof}
Let $F\subset \R^n$ be the facet into which $V$ is mapped,
and let $f:\R^n \to \R$ be the linear functional defining $F$, that is, such
that $f(x)=1$ for $x\in F$.
So, if $z_0 z_1\cdots$ is a word with letters in $V$,
then $f\after\phi(\overline z_0\cdots \overline z_{n-1})=n$ for all $n\in\Z$.
Let $y_0 y_1\cdots$ be a word with letters in $S$ such that
$\overline y_0\cdots \overline y_{m-1} = \overline z_0\cdots \overline z_{n-1}$
for some $n$ and $m$ in $\Z$. Since $f(x)\le 1$ for $x\in P$,
we have $m\ge f\after\phi(\overline y_0\cdots \overline y_{m-1}) = n$.
It follows that $z_0 z_1\cdots$ is a geodesic word.
\end{proof}

\begin{lemma}
\label{lem:rearrange}
Let $H$ be a nilpotent group generated as a group by a finite set $V$.
Then, for any $g\in H$, we can write $g\overline x = \overline y$,
where $x$ and $y$ are words with letters in $V$.
\end{lemma}
\begin{proof}
We will show the result is true for each of the groups $H/\Gamma_i$; $i\ge 2$,
using induction.

It is clearly true for the group $H/\Gamma_2= H/[H,H]$, since this group is
abelian.

Now assume the result is true for $H/\Gamma_i$ for some $i\ge 2$.

Let $C_j$ denote the set of simple commutators of weight $j$ of elements of $V$.
In other words, $C_1 := V$ and
\begin{align*}
C_{j+1} := \big\{ [c,s] \mid \text{$c\in C_j$ and $s\in V$} \big\},
\qquad \text{for $j\ge 1$}.
\end{align*}
From the induction hypothesis,
for every $c\in C_{i-1}$, there exists words $x_c$ and $y_c$ with letters in
$V$ such that $c\overline x_c  = \overline y_c $ mod $\Gamma_i$.
For each $c\in C_{i-1}$, let $g_c \in \Gamma_i$ be such that
\begin{align}
\label{eqn:interchange}
c\overline x_c  = \overline y_c  g_c.
\end{align}

Consider the word
\begin{align*}
w:= \prod_{\substack{c\in C_{i-1} \\ s\in V}}
                 (y_c s^{n(c,s)}x_c)(x_c s^{m(c,s)} y_c),
\end{align*}
where the product is taken in any fixed order, and $n(c,s)$ and $m(c,s)$
are non-negative integers for all $c\in C_{i-1}$ and $s\in V$.
Clearly, the letters of $w$ are in~$V$.

Using~(\ref{eqn:interchange}) and the fact that $[c,h]$ and $g_c$ are
central in $H/\Gamma_{i+1}$ for every $c\in C_{i-1}$ and $h\in H$, we get
\begin{align*}
\overline y_c  s^{n(c,s)}\overline x_c 
   &= c\overline x_c  g_c^{-1} s^{n(c,s)} c^{-1} \overline y_c  g_c,
          \qquad \mod \Gamma_{i+1} \\
   &= \overline x_c  [c,\overline x_c ] s^{n(c,s)} [c,s]^{n(c,s)}\overline y_c ,
          \qquad \mod \Gamma_{i+1}.
\end{align*}
Similarly,
\begin{align*}
\overline x_c  s^{m(c,s)}\overline y_c 
   &= \overline y_c  [\overline x_c ,c] s^{m(c,s)} [s,c]^{m(c,s)}\overline x_c ,
          \qquad \mod \Gamma_{i+1}.
\end{align*}
So
\begin{align*}
\overline w:= \prod_{\substack{c\in C_{i-1} \\ s\in V}}
                [c,s]^{n(c,s)-m(c,s)}
            \prod_{\substack{c\in C_{i-1} \\ s\in V}}
                (\overline x_c  s^{n(c,s)} \overline y_c)
                (\overline y_c  s^{m(c,s)} \overline x_c),
   \qquad \mod \Gamma_{i+1}.
\end{align*}
Since the group $\Gamma_i/\Gamma_{i+1}$ is generated by $C_i\mod\Gamma_{i+1}$,
and is in the center of $H/\Gamma_{i+1}$, any $h\in\Gamma_i$ can be 
written ($\!\!\!\!\mod \Gamma_{i+1}$) as the first factor in the
right-hand-side of the above equation,
with an appropriate choice of $n$ and $m$.
The second factor can be represented by a word $z$ with letters in $V$.
So $h\overline z = \overline w \mod \Gamma_{i+1}$.

Now let $g\in H$. By the induction hypothesis, $g\overline x = \overline y
\mod \Gamma_i$ for some words $x$ and $y$ with letters in $V$.
In other words, $g\overline x h=\overline y$ for some $h\in\Gamma_i$.
As we have just seen, $h\overline z = \overline w \mod \Gamma_{i+1}$
for some words $z$ and $w$ with letters in $V$.
So $g\overline x\overline w = \overline y \overline z \mod \Gamma_{i+1}$.
This completes the induction step.
\end{proof}

We are now ready to construct geodesic paths leading to the Busemann points
lying in finite orbits.

Let $F$ be a facet of $P$, and let $V$ be the set of generators mapped into
$F$ by $\phi$. We denote by $\langle V\rangle$ the subgroup of $G$ generated
as a group by the elements of $V$.
Of course, $\langle V\rangle$ is also nilpotent.
So, by Lemma~\ref{lem:rearrange}, for each $g\in\langle V\rangle$, there exist
words $x_g$ and $y_g$ with letters in $V$ such that
$g\overline x_g = \overline y_g$.
Let $C_i$ be the set of simple commutators of weight $i$ of elements of $V$.
Take a finite word $w_F$ with letters in $V$ such that
$x_c$ and $y_c$ appear as subwords, for every $c\in\bigcup_i C_i$.
The infinite word $w_F w_F \cdots$ is geodesic, by Lemma~\ref{lem:geo_facet}.
We define $\Lambda_F:\N\to G$ to be the geodesic path starting at the identity,
having steps given by this word. Since $\Lambda_F$ is geodesic, it converges
to a Busemann point of the horoboundary, which we denote by $\xi_F$.

\begin{lemma}
\label{lem:premult}
Let $F$ be a facet of $P$, and let $V$ be the set of generators mapped into
$F$ by $\phi$. Then, for any $g\in\langle V \rangle$, there is an $n\in\N$
such that $g\evalword^n$ can be written as a product of elements of $V$.
\end{lemma}
\begin{proof}
Let $\Gamma_i(\langle V\rangle)$ denote the decending central series of
$\langle V\rangle$, and let $C_i$ denote the set of simple commutators
of weight $i$ of elements of $V$.

We consider the groups $\langle V\rangle/\Gamma_i(\langle V\rangle)$; $i\ge 1$
and use induction.

The case when $i=1$ is trivial.

Now assume that $g\evalword^n$ can be written (mod $\Gamma_i(\langle V\rangle)$)
in the desired form.
That is to say, $g\evalword^n = \overline v c$ for some word $v$
with letters in $V$ and some $c\in\Gamma_i(\langle V\rangle)$.
We know that $c$ can be written as a product of simple commutators
and their inverses:
\begin{align*}
c = c_1^{\epsilon_1} \cdots c_m^{\epsilon_m},
\qquad\text{mod $\Gamma_{i+1}(\langle V\rangle)$},
\end{align*}
where $m\in\N$, and $c_j\in C_i$ and $\epsilon_j=\pm 1$
for all $j\in\{1,\dots,m\}$.
Recall that $\Gamma_i(\langle V\rangle)/\Gamma_{i+1}(\langle V\rangle)$
is in the center of $\langle V\rangle/\Gamma_{i+1}(\langle V\rangle)$.
Suppose that $\epsilon_j=1$ for some $j$.
Recall that $\word$ has subwords $x_{c_j}$ and $y_{c_j}$ satisfying
$c_j \overline x_{c_j} = \overline y_{c_j}$.
So by replacing the subword $x_{c_j}$ in
$\word$ by $y_{c_j}$, we obtain a word with letters in $V$ that represents
$c_j \evalword$ mod $\Gamma_{i+1}(\langle V\rangle)$.
Similarly, for any $j$ for which $\epsilon_j = -1$, we replace $y_{c_j}$
by $x_{c_j}$ to obtain a word with letters in $V$ representing
$c_j^{-1} \evalword$ mod $\Gamma_{i+1}(\langle V\rangle)$.
We concatenate all these words to get a word $u$ with letters in $V$
representing $c \evalword^m$ mod $\Gamma_{i+1}(\langle V\rangle)$.
So $vu$ represents $g\evalword^{n+m}$ mod $\Gamma_{i+1}(\langle V\rangle)$.
This concludes the induction step.
\end{proof}

\begin{lemma}
\label{lem:stabilises}
Let $F$ be a facet of $P$, and let $V$ be the set of generators mapped into
$F$ by $\phi$.
Then, the horofunction $\xi_F$ is invariant under the subgroup
$\langle V \rangle$.
\end{lemma}
\begin{proof}
Let $g\in\langle V \rangle$ and $n\in\N$.
Since $g[g,\evalword^n] \in \langle V \rangle$, we can find,
by Lemma~\ref{lem:premult},
$m\in\N$ such that $g[g,\evalword^n]\evalword^m$ is represented
by some word $v$ with letters in $V$.

Similarly,
we can find $l \in \N$ such that $[\evalword^{n+m},g]g^{-1}\evalword^{l}$
is represented by a word $v'$ with letters in $V$.

Consider the path $\gamma$ starting at the identity with steps given by the
word $\word^n v v' \word\word\cdots$.
By Lemma~\ref{lem:geo_facet}, it is a geodesic.
It agrees with the geodesic path $\facetgeo_F$ up to time $n|\word|$.
But $\evalword^n\overline v = g\evalword^{n+m}$, and
so $\gamma$ has a point in common with the path $g\facetgeo_F$ at
time $|\word^n v|$. We also have that
$\evalword^n \overline v \overline v' = \evalword^{n+m+l}$.
Therefore, $\gamma$ agrees with $\facetgeo_F$ after time $|\word^n v v'|$.

We now use Proposition~\ref{prop:switch} to deduce that $\facetgeo_F$ and
$g \facetgeo_F$ converge to the same limit $\xi_F$.
\end{proof}

\begin{theorem}
\label{thm:stabiliser}
Let $F$ be a facet of $P$, and let $V$ be the set of generators mapped into
$F$ by $\phi$. Then, the stabiliser of $\xi_F$ is $\langle V \rangle$.
\end{theorem}
\begin{proof}
We have already seen in Lemma~\ref{lem:stabilises} that $\langle V \rangle$
stabilises $\xi_F$.

Let $g\in G$ be such that $g\xi_F=\xi_F$. Consider the geodesics $\facetgeo_F$
and $g\facetgeo_F$. Since they have the same limit $\xi_F$, there exists
a geodesic $\gamma$ having infinitely many points in common with each, by
Proposition~\ref{prop:switch}.
So there exist $t_0$, $t_1$, and $t_2$ in $\N$ with $t_0\le t_1\le t_2$,
such that $\gamma(t_0)=\facetgeo_F(s_0)$, $\gamma(t_1)=g\facetgeo_F(s_1)$,
and $\gamma(t_2)=\facetgeo_F(s_2)$ for some $s_0$, $s_1$, and $s_2$ in $\N$.
Let $f:\R^n \to \R$ be the linear functional defining the facet $F$,
so that $f(x)=1$ for $x\in F$ and $f(x)<1$ for $x\in P\backslash F$.
Since the steps of $\facetgeo_F$ lie in $V$, we have
\begin{align*}
f\after\phi(\facetgeo_F(s_2)) - f\after\phi(\facetgeo_F(s_0)) = s_2 - s_0.
\end{align*}
But since $\gamma$ is a geodesic, $s_2 - s_0 = t_2 - t_0$.
We conclude that
\begin{align*}
f\after\phi(\gamma(t_2)) - f\after\phi(\gamma(t_0)) = t_2 - t_0.
\end{align*}
Using the fact that $f\after\phi$ is a homomorphism and that
$f\after\phi(x)\le 1$ for $x\in S$, we get that every step $x$ of $\gamma$
between $t_0$ and $t_2$ satisfies $f\after\phi(x)=1$, that is, they all
lie in~$V$.

Consider the finite path that starts at the identity, is identical to
$\facetgeo_F$ up to time $s_0$, follows $\gamma$ from
$\gamma(t_0)=\facetgeo_F(s_0)$ to $\gamma(t_1)=g\facetgeo_F(s_1)$, and
then follows $g\facetgeo_F$ backwards along to $g$. The steps of this path
lie entirely within $V\union V^{-1}$. We deduce that $g\in\langle V \rangle$.
\end{proof}

\begin{corollary}
For every facet $F$ of $P$, the Busemann point $\xi_F$ lies in a finite orbit.
\end{corollary}
\begin{proof}
Recall that $Z_0$ is the torsion subgroup of the abelianisation of $G$.
Choose a subset $R$ of $G$ containing one representative of each
equivalence class in $Z_0$.

Since $F$ is a facet of $P$, the set of vectors $\phi(V)$ generates a finite
index subgroup of $\Z^N$. So there is a finite subset $R'$ of $G$ such that,
for any $g\in G$, there exists $r'\in R'$ and $v\in\langle V\rangle$ satisfying
$\phi(g)=\phi(v)\phi(r')$.

So $\Theta(g) = \Theta(vr') \mod Z_0$, where $\Theta$ is the natural
homomorphism from $G$ to its abelianisation $G/[G,G]$.
This is equivalent to $\Theta(g) = \Theta(v)\Theta(rr')$ for some $r\in R$.
We conclude from this that $\Theta(\langle V\rangle)$ has finite index in
$\Theta(G)$. So, by Lemma~\ref{lem:group_finite_index}, $\langle V\rangle$
has finite index in $G$. But, by the theorem, $\langle V\rangle$ is exactly
the stabiliser of $\xi_F$, and so this Busemann point has a finite orbit.
\end{proof}

We have just shown that to every facet of $P$ there corresponds a finite
orbit of Busemann points. To prove the correspondence in the opposite
direction, we will need the following lemma. Note that its hypothesis that
the orbit is finite is not needed in the abelian setting~\cite{develin_cayley}.

\begin{lemma}
\label{lem:projects_to_facet}
Suppose the limit of a geodesic path $\gamma$ lies in a finite orbit.
Then, the projection map $\phi$ maps the direction of $\gamma$ into some
facet of $P$.
\end{lemma}
\begin{proof}
Let $u_0 u_1 \cdots$ be the geodesic word giving the steps of $\gamma$.
For each $n\in\N$, let $\gamma^n$ be the geodesic starting at the identity
and having steps given by the word $u_n u_{n+1} \cdots$.
The limit of each $\gamma^n$ is in the same orbit as that of $\gamma$,
which we have assumed finite. Therefore, there is some subsequence
of geodesics $\gamma^{n_i}$ all of whose limits are equal.
Call the common limit $\xi$.

For each $i\in\N$, define the finite word $w_i:=u_{n_i} \cdots u_{n_{i+1}-1}$,
and let $W:=\{w_i \mid i\in\N \}$ be the set of all these words.
Fix $i\in\N$. Since $\gamma^{n_i}$ is eventually just a copy of
$\overline w_i \gamma^{n_{i+1}}$ shifted in time, these two geodesics have
the same limit. But both $\gamma^{n_i}$ and $\gamma^{n_{i+1}}$ converge
to $\xi$. Therefore, $w_i$ represents a group element in the stabiliser
of $\xi$.

Now consider an infinite concatenation $z := z_0 z_1 \cdots z_m \cdots$ of
words $z_m$; $m\in\N$, each in $W$. Given any $m\in\N$, the word $z$ can be
written
\begin{align*}
z = z_0 \cdots z_{m-1} u_{n_i} \cdots u_{n_{i+1}-1} z_{m+1} \cdots
\end{align*}
because $z_m = w_i$ for some $i\in\N$.
Since $p:= \overline z_0 \cdots \overline z_{m-1}$
is in the stabiliser of $\xi$,
we have $\xi(pg) = \xi(p) + \xi(g)$ for all $g\in G$.
So
\begin{align*} 
\xi(\overline z_0  \cdots \overline z_{m-1} \overline u_{n_i}
          \cdots \overline u_{n_i +l-1})
   = \xi(\overline z_0 \cdots \overline z_{m-1})
        + \xi(\overline u_{n_i} \cdots \overline u_{n_i +l-1}),
\end{align*}
for $l\in\{0,\ldots,n_{i+1}-n_i\}$.
But since $\gamma^{n_i}$ is a geodesic converging to $\xi$, we have
$\xi(\overline u_{n_i} \cdots \overline u_{n_i +l-1})= -l$.
Since this reasoning is valid for any $m\in\N$, we conclude that $\xi$,
evaluated along the path $\lambda$, decreases at every step, where
$\lambda$ is the path starting at the identity with steps given by the word
$z$. It follows from this that $\lambda$ is a geodesic.
We have thus shown that every concatenation of words from $W$ is a geodesic
word.

Take now sufficiently many words $z_0, \ldots, z_{n-1}$ from $W$ that every
generator in the direction of $\gamma$ appears in at least one.
Let $a:=\phi(\overline z_0 \cdots \overline z_{n-1} )$,
and let $b$ be the intersection
of the line $0a$ in $\R^N$ with the boundary of $P$. Of course, $b$ has
rational coordinates and lies in some facet $F$ of $P$.
So, by Lemma~7 of~\cite{develin_cayley}, $b$ may be expressed as an affine
combination of the vertices of $F$, with non-negative rational coefficients.
Since $a$ is a rational multiple of $b$, we may express $a$ in the same way.
Multiplying through by the common denominator, we can find a positive integer
$l$ as large as we wish, such that $la$ is a linear combination of the vertices
of $F$ with non-negative integer coefficients. So there exists a word $x$ with
letters in $\phi^{-1}(F) \intersection S$ such that $\phi(\overline x)=la$.

Recall that $Z_0$ is the torsion subgroup of the abelianisation $G/[G,G]$.
Let $T\subset G$ be such that there is exactly one element of $T$ in each
of the equivalence classes in $Z_0$. Of course, $|T| = |Z_0| < \infty$.

We can now write $(\overline z_0 \cdots \overline z_{n-1})^l = \overline xtc$,
with $t\in T$ and $c\in[G,G]$.

We have assumed that the stabiliser of $\xi$, which
we now call $U$, has finite index in $G$.
Therefore, by Lemma~\ref{lem:comm_finite_index}, $[U,U]$ has finite index in
$[G,G]$. So we can write $c=ru$, with $r$ in some finite subset $R$ of $G$,
and $u$ in $[U,U]$. Because $U$ stabilises $\xi$, the restriction of $\xi$
to $U$ is a homomorphism into $\Z$. In particular, $\xi(u^{-1})=0$.
So, for $i\in\N$ large enough,
$d(u^{-1},\gamma^{n_0}(i))= d(e,\gamma^{n_0}(i)) = i$.

Let
\begin{align*}
D:= \max\{ d(e,t'r') \mid \mbox{$t'\in T$ and $r' \in R$} \}.
\end{align*}
Of course, $D$ is finite.
We have
\begin{align*}
d(e,\overline x tru \gamma^{n_0}(i))
         & \le d(e,\overline x) + d(e,tr) + d(e,u\gamma^{n_0}(i)) \\
         & \le |x| + D + i.
\end{align*}
But $|x| = l f(a)$, where $f:\R^N\to \R$ is the linear functional such that
$f(g)=1$ for all $g\in F$.
Also, since $(z_0 \cdots z_{n-1})^l w_0 w_1 w_2 \cdots$ is a geodesic word,
\begin{align*}
d(e,(\overline z_0 \cdots \overline z_{n-1})^l\gamma^{n_0}(i))
                   = l|z_0 \cdots z_{n-1}| + i.
\end{align*}
Therefore  $l f(a) + D \ge l |z_0 \cdots z_{n-1}|$.
Since $l$ can be made arbitrarily large, $f(a) \ge |z_0 \cdots z_{n-1}|$.
We may now apply the fact that $f\after \phi(s) \le 1$ for every $s\in S$
to deduce that $f\after \phi(g)= 1$ for every letter $g$ in
$z_0 \cdots z_{n-1}$, which includes every letter in the direction of $\gamma$.
Therefore, all these letters are mapped by $\phi$ to $F$.
\end{proof}

\begin{theorem}
\label{all_finite_orbits}
Every finite orbit of Busemann points can be written
$\{g\xi_F \mid g\in G\}$ for some facet $F$ of $P$.
\end{theorem}
\begin{proof}
Let $\gamma_1$ be a geodesic converging to a point in a finite orbit.
By Lemma~\ref{lem:projects_to_facet}, the image of the direction of $\gamma_1$
under $\phi$ lies in some facet $F$ of $P$. By, if necessary, premultiplying
by an appropriate factor and removing an initial section of the path, we can
find a geodesic path $\gamma$ converging to a point $\xi$ in the same orbit
as the limit of $\gamma_1$ that starts at the identity and for which the image
under $\phi$ of every step lies in $F$.
Let the word $\word$ be defined as above. Since $\xi$ is in a finite orbit,
there exists some $m\in\N$ with $m>0$ such that $\evalword^m\xi$ is equal
to $\xi$. So $\evalword^{mn}\xi=\xi$ for all $n\in\N$.
Let $z_0 z_1 \cdots$ be the geodesic word giving the steps of $\gamma$,
and let $\gamma_2$ be the geodesic path starting at the identity and having
steps given by $\word^{mn}z_0 z_1 \cdots$ for some large $n\in\N$.
By Lemma~\ref{lem:geo_facet}, $\gamma_2$ is a geodesic. Also, $\gamma_2$ agrees
with $\facetgeo_F$ up to time $nm|\word|$, which can be made as large as one
likes. Finally, $\gamma_2$ converges to $\xi$, and so has the same limit as
$\gamma$.

We have proved one half of condition~(iv) of Proposition~\ref{prop:switch}.
To prove the other, take $n$ as large as you want and let
$g:= \overline z_0 \cdots \overline z_n$.
Clearly, $g^{-1} \in \langle V \rangle$,
where $V$ is the set of generators mapped into $F$ by $\phi$.
Therefore, by Lemma~\ref{lem:premult}, there is some $m\in\N$ such that
$g^{-1} \evalword^m$ can be represented by some word $v$ with letters in $V$.
By Lemma~\ref{lem:geo_facet}, the path starting at the identity and having
steps given by $z_0 \cdots z_n v \word \word \cdots$ is a geodesic.
This path agrees with $\gamma$ up to time $n$, and it eventually agrees with
$\facetgeo_F$ because $\overline z_0 \cdots \overline z_n \overline v = \evalword^m$.

It now follows from Proposition~\ref{prop:switch} that $\gamma$ and
$\facetgeo_F$ have the same limit, in other words, that $\xi=\xi_F$.
\end{proof}

We finish this section by showing how the correspondence in the abelian setting
between proper faces of $P$ and orbits of Busemann points breaks down for
general finitely-generated nilpotent groups.

\textbf{Example 1.}
Consider the following group of nilpotency class 3:
\begin{align*}
G:= \langle a,b \mid [a,g] = [b,g] = [a,h] = [b,h] = e \rangle,
\end{align*}
where $c=[a,b]$, $g:=[a,c]$, and $h:=[b,c]$.

Any element of $G$ can be written in the form $g^i h^j c^k b^l a^m$,
with $i$, $j$, $k$, $l$, and $m$ in $\Z$. We may take the projection map $\phi$
to be $\phi(g^i h^j c^k b^l a^m) := (l,m)$. The image of the generating set
$\{a,b,a^{-1},b^{-1}\}$ is a square. Let $F$ be the facet defined by the
corners $\{(1,0), (0,1)\}$.

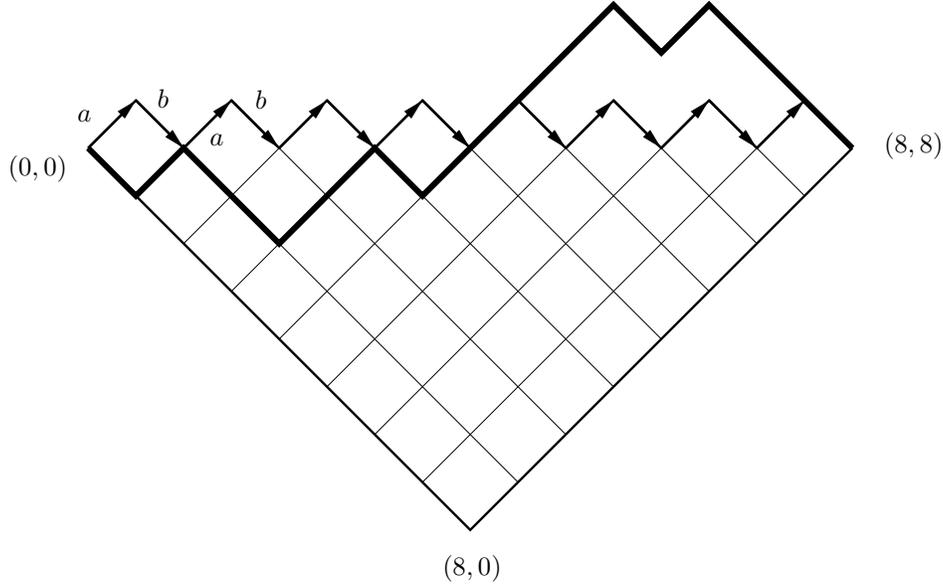
\begin{figure}
\input{water.pstex_t}
\caption{Diagram for the proof of Lemma~\ref{lem:example_perm}.}
\label{fig:water}
\end{figure}

Define the word $x:= ab$. We will show that, although $xxx\cdots$ and
$\facetgeo_F$ have the same directions, their limits differ.
\begin{lemma}
\label{lem:example_perm}
Let $w$ be any word with letters in $\{a,b\}$ such that
$\overline w = g^m h^n {\overline x}^l$ for some $m,n \in \Z$ and $l\in\N$.
Then, $m-n$ is positive unless $w=x^l$.
\end{lemma}
\begin{proof}
Represent $w$ as a piecewise affine curve that starts at the origin and takes
a unit step in the $x$-direction for each letter $b$, and unit step in the
$y$-direction for each letter $a$.
The curve will obviously finish at $(l,l)$.
Write $\overline w$ in the form $g^i h^j c^k b^l a^l$.
It is not hard to see that $k$ is the number of unit squares under the curve.
After further thought, one also sees that $j$ is the sum of the $x$-coordinates
of the upper right corners of these squares, and that $i$ is the sum of the
$y$-coordinates.
So the center of gravity, taking the mass of each unit square to be
concentrated on its upper right corner, is $(j/k, i/k)$.
We now redraw the figure with the axes rotated through a $45^\circ$ angle,
as in Figure~\ref{fig:water}. We see that choosing $w$ to minimise $i-j$
is the same as finding the curve between $(0,0)$ and $(l,l)$ having $(l+1)l/2$
unit squares underneath it that minimises the height of the center of mass
of these squares. But it is obvious that there is a unique minimising
curve, namely that obtained from $x^l$. The conclusion follows.
\end{proof}

\begin{prop}
\label{prop:example}
The limit $\eta$ of $xxx\cdots$ is different from $\xi_F$.
\end{prop}
\begin{proof}
Let $l\in\N$ and let $w$ be a geodesic word with letters in
$\{a,b,a^{-1},b^{-1}\}$ representing $g^{-1} \overline x^l$.
By Lemma~\ref{lem:example_perm}, either $a^{-1}$ or $b^{-1}$ occurs as a
letter of $w$. Also, the excess of $a$s over $a^{-1}$s in $w$ equals $l$,
as does the excess of $b$s over $b^{-1}$s.
Therefore, $d(e, g^{-1}\overline x^l) = |w| > 2l$.
We conclude that
\begin{align*}
\eta(g) = \lim_{l\to\infty} \Big( d(g,\overline x^l) - d(e,\overline x^l) \Big)
        >0.
\end{align*}

However, $\{a,b\}$ generates $G$, and so $\xi_F$ is a fixed point
by Theorem~\ref{thm:stabiliser}. So $\xi_F$ is a homomorphism from $G$ to $\Z$.
In particular, $\xi_F(g) = 0$. Therefore, $\eta$ and $\xi_F$ differ.
\end{proof}

\section{The boundary of the discrete Heisenberg group}

The discrete Heisenberg group $H_3$ is the group of $3\times 3$ upper triangular
matrices of the form 
\begin{align*}
\begin{pmatrix} 1 & x & z\\ 0 & 1 & y\\ 0 & 0 & 1\\ \end{pmatrix},
\qquad\text{with $x$, $y$, and $z$ in $\Z$}.
\end{align*}
It is the simplest non-abelian nilpotent group. Let
\begin{align*}
a:=\begin{pmatrix} 1 & 1 & 0\\ 0 & 1 & 0\\ 0 & 0 & 1\\ \end{pmatrix},
\qquad
b:=\begin{pmatrix} 1 & 0 & 0\\ 0 & 1 & 1\\ 0 & 0 & 1\\ \end{pmatrix},
\qquad\text{and}\quad
c:=\begin{pmatrix} 1 & 0 & 1\\ 0 & 1 & 0\\ 0 & 0 & 1\\ \end{pmatrix}.
\end{align*}
Observe that $ab=bac$, that is, $c$ is the commutator of $a$ and $b$,
and that the center of $H_3$ is the cyclic group generated by $c$.
The following is a presentation of $H_3$:
\begin{align*}
H_3 := \langle a,b \mid [[a,b],a] = [[a,b],b] = e \rangle.
\end{align*}
All elements of $H_3$ can be written in the form $c^z b^y a^x$,
with $x,y,z\in\Z$.

We will calculate all the Busemann points of this group with the word length
metric coming from the standard generating set $S:=\{a,b,a^{-1},b^{-1}\}$.
Our work relies on the following formula for this metric,
found by Blach{\`e}re~\cite{blachere}.
\begin{theorem}[\cite{blachere}]
\label{thm:formula}
Let $g:= c^z b^y a^x$ be an element of $H_3$.
Its distance to the identity with respect to the generating set
$S:=\{a,b,a^{-1},b^{-1}\}$ is
\begin{itemize}
\renewcommand{\labelitemi}{I.}
\item
if $zyx\ge 0$, and
   \begin{itemize}
   \renewcommand{\labelitemii}{I.1.}
   \item
   if $\max(x^2,y^2) \le |z|$, then
      \begin{align*}
      d(e,g) = 2\big\lceil 2\sqrt{|z|}\big\rceil - |x| - |y|;
      \end{align*}
   \renewcommand{\labelitemii}{I.2.}
   \item
   if $\max(x^2,y^2) \ge |z|$, and
       \begin{itemize}
          \renewcommand{\labelitemiii}{I.2.1.}
          \item
          $|xy| \ge |z|$, then
            \begin{align*}
            d(e,g) = |x| + |y|;
            \end{align*}
          \renewcommand{\labelitemiii}{I.2.2.}
          \item
          $|xy| \le |z|$, then
            \begin{align*}
            d(e,g) = 2\lceil \min(|z/x|,|z/y|) \rceil + \big| |x| - |y| \big|;
            \end{align*}
       \end{itemize}
   \end{itemize}
\renewcommand{\labelitemi}{II.}
\item
if $zyx\le 0$, and
   \begin{itemize}
   \renewcommand{\labelitemii}{II.1.}
   \item
   if $\max(x^2,y^2) \le |z| + |xy|$, then
      \begin{align*}
      d(e,g) = 2\big\lceil 2\sqrt{|z|+|xy|}\big\rceil - |x| - |y|;
      \end{align*}
   \renewcommand{\labelitemii}{II.2.}
   \item
   if $\max(x^2,y^2) \ge |z| + |xy|$, then
         \begin{align*}
         d(e,g) = 2\lceil \min(|z/x|,|z/y|) \rceil + |x| + |y|.
         \end{align*}
   \end{itemize}
\end{itemize}
\end{theorem}

Given a word $y$, we  call the pairs of consecutive letters of $y$
the \emph{transitions} of $y$. The transitions $ab$, $ba^{-1}$, $a^{-1}b^{-1}$,
and $b^{-1}a$ are said to be positive. The transitions obtained by reversing
the letters of these are called negative.
Observe that the group element corresponding to a positive transition is equal
to $c$ times the element corresponding to the negative transition obtained
by reversing the order of the letters. So, if one takes a word and transforms
it by reversing an equal number of positive and negative transitions,
the resulting word is a geodesic if and only if the original was.

\begin{prop}
\label{prop:h3geodesics}
Every infinite geodesic word in $H_3$ either has two elements of $S$
as letters or consists of a finite prefix followed by a single letter repeated.
\end{prop}
\begin{proof}
Let $y$ be an infinite word that has at least three elements of $S$ as letters,
at least two of which appear infinitely often.

We wish to show that $y$ is not a geodesic.
Obviously, this is the case if it is not freely reduced, so we assume
the contrary. Our above assumptions on $y$ then imply that there are either
infinitely many positive transitions in $y$ or there are
infinitely many negative transitions.

Since $y$ has at least three elements of $S$ as letters, it must have
as letters both some generator and its inverse.
By repeatedly swapping this generator with one of its neighbours,
we can bring it to a position adjacent to its inverse.
If $y$ has both infinitely many positive and infinitely many negative
transitions,
then we can compensate for these swaps by reversing transitions elsewhere
so that the total number of positive transitions reversed equals the
total number of negative ones reversed.
Since the resulting word is clearly not a geodesic,
it follows that $y$ is not either.

So we will assume without loss of generality that there are infinitely
many positive transitions and finitely many negative ones.
Therefore $y$ can be written in the form
\begin{align*}
y = R a^{i_0} b^{j_0} a^{-i_1} b^{-j_1} a^{i_2}\cdots,
\end{align*}
where $(i_n)$ and $(j_n)$ are sequences of positive integers, and $R$ is some
finite initial word.

Consider what happens if $(i_n)$ and $(j_n)$ eventually become, respectively,
constants $i$ and $j$. Define the word $w:=a^i b^j a^{-i} b^{-j}$.
Observe that $\overline w = c^{ij}$. So, the word $www\cdots$ is not geodesic
since, by the formula of Theorem~\ref{thm:formula}, $d(e,\overline w^n)$
grows like the square root of $n$, rather than linearly.
Therefore, in this case, $y$ is not geodesic.

If $(i_n)$ and $(j_n)$ are not both eventually constant, then
we can find $n\ge 1$ such that $i_n < i_{n+1}$ or $j_n < j_{n+1}$.
Assume that the latter holds; the former case can be handled similarly.
Also assume without loss of generality that $n$ is even.
So $y$ can be written
\begin{align*}
y = R \cdots b^{-j_{n-1}} (a^{i_n}) b^{j_n} a^{-i_{n+1}} b^{-j_{n+1}} \cdots.
\end{align*}
By moving the bracketed clump of $a$'s forward, we get a word
\begin{align*}
z = R \cdots b^{-j_{n-1}} b^{j_n} a^{-i_{n+1}}
                             b^{-j_n} (a^{i_n}) b^{-j_{n+1}+j_n} \cdots,
\end{align*}
which is clearly not a geodesic.
However, $z$ is obtained from $y$ by reversing the letters of $i_n j_n$
positive transitions and then reversing the letters of the same number of
negative transitions. It follows that $y$ is not geodesic either.
\end{proof}

\newcommand\epsilona{{\epsilon_a}}
\newcommand\epsilonb{{\epsilon_b}}
\newcommand\sign{\operatorname{sign}}

\newcommand\corner{\xi}
\newcommand\ahoro{\eta}
\newcommand\bhoro{\zeta}

\begin{theorem}
\label{thm:corners}
Let $\epsilona$ and $\epsilonb$ be in $\{-1,+1\}$, and let $\Omega$ be the set
of all words having letters in $\{a^\epsilona,b^\epsilonb\}$ with both letters
occurring infinitely often.
Then, all words in $\Omega$ are geodesic, and the corresponding paths
all converge to the same Busemann point, given by
\begin{align*}
\corner^{\epsilona \epsilonb}(c^z b^y a^x) := -\epsilona x - \epsilonb y.
\end{align*}
\end{theorem}
\begin{proof}
Let $w$ be a finite word having set of letters $\{a^\epsilona,b^\epsilonb\}$,
and let $\overline w=c^z b^y a^x$ be the associated element of the group.
If $z\neq 0$, then $\sign z = \sign y. \sign x = \epsilonb \epsilona$.
Therefore, $zyx\ge0$.
Also, $|z|\le |xy|$. We deduce from the word length formula of
Theorem~\ref{thm:formula} that $d(e,\overline w) = |x| + |y|$, which is exactly
the number of letters in $w$. So $w$ is geodesic. It follows that every
word in $\Omega$ is geodesic.

Let $w_1$ and $w_2$ be words in $\Omega$, and let $\gamma_1$ and $\gamma_2$
be the corresponding geodesic paths starting at the identity.
Choose $N\in\N$, and let $u$ be the prefix of $w_1$ of length $N$.
Since $w_2$ has each of the letters $a^\epsilona$ and
$b^\epsilonb$ infinitely often, we can find a prefix $v$ of $w_2$ that
has more $a^\epsilona$'s and $b^\epsilonb$'s than $u$.
So there exists a word $x$ with letters in $\{a^\epsilona,b^\epsilonb\}$
of length $|v|-|u|$ such that $\overline u \overline x = c^m \overline v$
for some $m\in\Z$. Since $w_2$ has infinitely many positive transitions
and infinitely many negative transitions, we can write it $w_2=vy\omega$,
where $\omega$ is an infinite word (in $\Omega$) and $y$ is a finite word
satisfying $\overline y = c^m \overline z$, for some $z$ with letters
in $\{a^\epsilona,b^\epsilonb\}$.
By the first part of the proposition, $uxz\omega$ is a geodesic word.
The corresponding geodesic path agrees with $\gamma_1$ up to
time $N$ and agrees with $\gamma_2$ after time $|uxz|$. A similar argument
shows that there is also a geodesic path agreeing with $\gamma_2$ up to time
$N$ and eventually coinciding with $\gamma_1$.
We now apply Proposition~\ref{prop:switch} to deduce that $\gamma_1$
and $\gamma_2$ converge to the same limit.

Consider the sequence $h_n := (a^\epsilona b^\epsilonb)^n$.
By the previous part of the proposition, this sequence converges to the common
limit of the geodesics in $\Omega$. We have
\begin{align*}
h_n = c^{\epsilona\epsilonb n(n+1)/2} b^{\epsilonb n} a^{\epsilona n}.
\end{align*}
Let $g:= c^z b^y a^x$ be an arbitrary point in $H_3$. We have
\begin{align*}
d(g,h_n) = d(e,g^{-1} h_n)
      = d(e,c^{yx-z-\epsilonb xn + \epsilona\epsilonb n(n+1)/2}
            b^{\epsilonb n - y}
            a^{\epsilona n - x}).
\end{align*}
Consider what happens when $n$ is large. The exponents of $c$, $b$, and $a$
will have the signs $\epsilona\epsilonb$, $\epsilonb$, and $\epsilona$,
respectively; therefore, their product will be positive.
The square of the exponents of both $a$ and $b$ will be approximately
twice the absolute value of that of $c$, as will the absolute value of
the product of the exponents of $a$ and $b$.
Looking at the formula of Theorem~\ref{thm:formula},
we see that, for large $n$, the relevant case is I.2.1,
and so $d(g,h_n)=|\epsilonb n - y| + |\epsilona n - x|$.
In particular $d(e,h_n)=2n$. So the limiting horofunction is
$\corner^{\epsilona \epsilonb}(c^z b^y a^x) := -\epsilona x - \epsilonb y$.
\end{proof}

\begin{theorem}
\label{thm:heisenberg_busemen}
The following functions are Busemann points of $H_3$:
\begin{align*}
\ahoro^{\epsilon}_{m,n}(c^k b^j a^i)
   &:= -\epsilon i + |j-n| - |n|
        + 2 J \big(\epsilon (j - n), (j-n)i - (k-m)\big) - 2J(-\epsilon n,m), \\
\bhoro^{\epsilon}_{m,l}(c^k b^j a^i)
   &:= -\epsilon j + |i-l| - |l|
        + 2 J \big(-\epsilon (i - l), jl - (k-m)\big) - 2J(\epsilon l,m),
\end{align*}
where $\epsilon\in\{-1,+1\}$, $m, n, l\in\Z$, and
\begin{align*}
J(u,v) := \begin{cases}
   1, & \text{if $v\neq 0$ and $uv\ge 0$,} \\
   0, & \text{otherwise}, \\
   \end{cases}
\qquad\text{for $u,v\in\Z$}.
\end{align*}
Together with $\corner^{++}$, $\corner^{-+}$, $\corner^{++}$,
and $\corner^{--}$, these are the only Busemann points.
\end{theorem}
\begin{proof}
Choose $\epsilon\in\{-1,+1\}$ and $m, n\in\Z$, and define the path
$\gamma^\epsilon_{m,n}(t) := c^m b^n a^{\epsilon t}$, for $t\in\Z$.
By the formula of Theorem~\ref{thm:formula},
\begin{align*}
d(\gamma^\epsilon_{m,n}(0),\gamma^\epsilon_{m,n}(t)) = d(e,a^{\epsilon t}) = t.
\end{align*}
Therefore $\gamma^\epsilon_{m,n}$ is a geodesic starting at $c^m b^n$.

Let $g:=c^k b^j a^i$ be an element of $H_3$. We have
\begin{align*}
d(g,\gamma^\epsilon_{m,n}(t)) = d(e,c^{m-k+(j-n)i} b^{n-j} a^{\epsilon t-i}).
\end{align*}
Write
\begin{align*}
x_t:=\epsilon t-i, \qquad
y:= n-j, \qquad
\text{and $z:= m-k + (j-n) i$}.
\end{align*}
So $|x_t|$ tends to infinity as $t$ tends to infinity, whereas $y$ and $z$
remain constant.

Therefore the only relevant cases in the formula of Theorem~\ref{thm:formula}
when $t$ is large are I.2.1.~and II.2.
If $z=0$, then both cases give the same result:
$|x_t|+|y|$. If $z\neq 0$, then we have case I.2.1~if $\epsilon yz > 0$,
in which case $d(g,\gamma^\epsilon_{m,n}(t))=|x_t|+|y|$,
or case II.2.~if $\epsilon yz\le 0$, which
gives $d(g,\gamma^\epsilon_{m,n}(t))=|x_t|+|y|+2$, for $x_t$ large enough.
To sum up:
\begin{align*}
d(g,\gamma^\epsilon_{m,n}(t)) = |x_t|+|y|+2 J(-\epsilon y,z).
\end{align*}
Subtracting from this the value obtained when $i=j=k=0$, and taking the limit
as $t$ tends to infinity, we see that the path $\gamma^\epsilon_{m,n}$ converges
to the horofunction $\ahoro^{\epsilon}_{m,n}$ defined above.
This horofunction is a Busemann point since $\gamma^\epsilon_{m,n}$ is a
geodesic.

Similar reasoning shows that the path
$\lambda^\epsilon_{m,l}(t):= c^{m+\epsilon tl} b^{\epsilon t} a^l$ converges to
$\bhoro^\epsilon_{m,l}$, for all $\epsilon\in\{-1,+1\}$ and  $m, l\in\Z$, and
that the limit is a Busemann point.

By Proposition~\ref{prop:h3geodesics}, every infinite geodesic in $H_3$ either
uses exactly two generators or uses all but one only a finite number
of times. In the former case, the geodesic converges to one of 
$\corner^{++}$, $\corner^{-+}$, $\corner^{+-}$, and $\corner^{--}$,
by Theorem~\ref{thm:corners}.
In the latter, the geodesic must eventually coincide with either
$\gamma^\epsilon_{m,n}$ or $\lambda^\epsilon_{m,l}$ for some value of
$\epsilon$ and of $m$ and $n$, or $m$ and $l$. Therefore, in this case,
the geodesic converges to the corresponding $\ahoro^{\epsilon}_{m,n}$
or $\bhoro^\epsilon_{m,l}$.
\end{proof}

The action of $H_3$ on each of the Busemann points can easily be calculated:
\begin{align*}
(c^z b^y a^x) \ahoro^\epsilon_{m,n} &= \ahoro^\epsilon_{m+z+nx,n+y}, \\
(c^z b^y a^x) \bhoro^\epsilon_{m,l} &= \bhoro^\epsilon_{m+z-y(l+x),l+x}, \\
(c^z b^y a^x) \corner^{\epsilona,\epsilonb} &= \corner^{\epsilona,\epsilonb}.
\end{align*}
So, in particular, each of the Busemann points $\corner^{++}$, $\corner^{+-}$,
$\corner^{-+}$, and $\corner^{--}$ is fixed by the action of $H_3$.
We have here an illustration of Theorem~\ref{thm:maintheorem}.
The abelianisation of $H_3$ is isomorphic to $\Z^2$. One can take the map
$\phi$ to be $\phi(c^z b^y a^x) := (x,y)$. The image of $S$ under this map is
the set $\{(1,0),(0,1),(-1,0),(0,-1)\}$, and so $P$ is the square with these
points as corners. Of course, $P$ has four facets, corresponding to the four
fixed Busemann points.

In \cite{winweb_busemann}, Webster and Winchester conjectured that 
there is a boundary point or points of the form
$\lim_{i\to\infty}\ahoro^+_{m_i,n_i}$,
where $m_i$ and $n_i$ grow without bound as $i$ tends to infinity,
with $m_i\ge \alpha n_i$ eventually for any $\alpha$, and that this point
or these points are fixed under the action of $H_3$.
Using the formula for $\ahoro^+_{m,n}$ from Theorem~\ref{thm:heisenberg_busemen}
one can calculate that there is exactly one point of this form,
namely $\corner^{++}$, which is indeed a fixed point.

\bibliographystyle{plain}
\bibliography{nilpotent}

\end{document}

%% file: water.pstex_t
\begin{picture}(0,0)%
\includegraphics{water.pstex}%
\end{picture}%
\setlength{\unitlength}{1973sp}%
\begingroup\makeatletter\ifx\SetFigFont\undefined%
\gdef\SetFigFont#1#2#3#4#5{%
  \reset@font\fontsize{#1}{#2pt}%
  \fontfamily{#3}\fontseries{#4}\fontshape{#5}%
  \selectfont}%
\fi\endgroup%
\begin{picture}(11840,7316)(796,-8212)
\put(2671,-2251){\makebox(0,0)[lb]{\smash{{\SetFigFont{10}{12.0}{\rmdefault}{\mddefault}{\updefault}{\color[rgb]{0,0,0}$b$}%
}}}}
\put(3331,-2731){\makebox(0,0)[lb]{\smash{{\SetFigFont{10}{12.0}{\rmdefault}{\mddefault}{\updefault}{\color[rgb]{0,0,0}$a$}%
}}}}
\put(1666,-2431){\makebox(0,0)[lb]{\smash{{\SetFigFont{10}{12.0}{\rmdefault}{\mddefault}{\updefault}{\color[rgb]{0,0,0}$a$}%
}}}}
\put(796,-3106){\makebox(0,0)[lb]{\smash{{\SetFigFont{10}{12.0}{\rmdefault}{\mddefault}{\updefault}{\color[rgb]{0,0,0}$(0,0)$}%
}}}}
\put(3901,-2266){\makebox(0,0)[lb]{\smash{{\SetFigFont{10}{12.0}{\rmdefault}{\mddefault}{\updefault}{\color[rgb]{0,0,0}$b$}%
}}}}
\put(11806,-2821){\makebox(0,0)[lb]{\smash{{\SetFigFont{10}{12.0}{\rmdefault}{\mddefault}{\updefault}{\color[rgb]{0,0,0}$(8,8)$}%
}}}}
\put(6256,-8131){\makebox(0,0)[lb]{\smash{{\SetFigFont{10}{12.0}{\rmdefault}{\mddefault}{\updefault}{\color[rgb]{0,0,0}$(8,0)$}%
}}}}
\end{picture}%